\documentclass[a4paper]{article}

\usepackage[numbers,sort&compress]{natbib}
\usepackage{amsmath,amssymb,amsfonts,amsthm}
\usepackage{graphicx}
\usepackage{xcolor}
\usepackage{microtype}
\usepackage{authblk}

\AtBeginDocument{%
	\setlength{\abovedisplayskip}{6pt plus 2pt minus 2pt}%
	\setlength{\belowdisplayskip}{6pt plus 2pt minus 2pt}%
	\setlength{\abovedisplayshortskip}{3pt plus 2pt minus 1pt}%
	\setlength{\belowdisplayshortskip}{4pt plus 2pt minus 2pt}%
	\setlength{\textfloatsep}{10pt plus 2pt minus 2pt}%
	\setlength{\floatsep}{8pt plus 2pt minus 2pt}%
	\setlength{\intextsep}{8pt plus 2pt minus 2pt}%
	\setlength{\jot}{2pt}%
}

\DeclareMathOperator{\rk}{rk}

\theoremstyle{plain}
\newtheorem{theorem}{Theorem}
\newtheorem{lemma}[theorem]{Lemma}
\newtheorem{proposition}[theorem]{Proposition}
\newtheorem{corollary}[theorem]{Corollary}

\theoremstyle{definition}
\newtheorem{definition}[theorem]{Definition}

\theoremstyle{remark}
\newtheorem*{remark}{Remark}

\title{Asymptotic of Betti numbers of the resonance arrangement}

\author[1]{Nikolai K. Zuev}
\author[1]{Mikhail V. Bludov\thanks{Corresponding author.}}
\affil[1]{Laboratory of Combinatorial and Geometric Structures, Moscow Institute of Physics and Technology, 9 Institutskiy per., Dolgoprudny, 141701, Moscow Region, Russian Federation}
\affil[ ]{\textit{E-mail addresses:} \texttt{zuev.nk@phystech.edu} (N.~K.~Zuev), \texttt{bludov.mv@phystech.edu} (M.~V.~Bludov)}

\date{}

\begin{document}
	
	\maketitle
	
	\begin{abstract}
		Betti numbers are classical topological invariants. Rather than considering general topological spaces, we focus on arrangements of hyperplanes---specifically, on the resonance arrangement $\mathcal{A}_n$, which in a precise sense is universal among rational arrangements. We study the asymptotic of $b_m(\mathcal{A}_n)$ for every fixed $m$ and variable $n$. Our first result is the estimate
		\[
		b_m(\mathcal{A}_n) = \frac{1}{m!}2^{nm} - O\!\left(\left(3 \cdot 2^{m - 2}\right)^n\right).
		\]
		Our second result refines it: the exponential order of the correction term is attained, and its coefficient is explicit,
		\[
		b_m(\mathcal{A}_n)
		=
		\frac{2^{nm}}{m!}
		-
		\frac{m+2}{8(m-2)!}
		\bigl(3\cdot2^{m-2}\bigr)^n
		+
		O_m\!\left(
		\bigl(3\cdot2^{m-2}-1\bigr)^n
		\right).
		\]
		The proof combines K\"uhne's theory of prototypes with a new matrix-based description of minimal balanced collections.
	\end{abstract}
	
	\noindent\textbf{Keywords:} betti numbers, resonance arrangement, column complementation, asymptotic enumeration, balanced collection, $0$--$1$ matrix
	
	\section{Introduction}\label{sec:intro}
	
	The \emph{resonance arrangement} $\mathcal{A}_n$ is the hyperplane arrangement in $\mathbb{R}^n$ whose hyperplanes are
	\[
	H_I=\Bigl\{x\in\mathbb{R}^n:\sum_{i\in I}x_i=0\Bigr\},
	\qquad \varnothing\neq I\subseteq[n].
	\]
	It is, in a precise sense, universal among rational hyperplane arrangements: every rational arrangement is a minor of $\mathcal{A}_n$ for sufficiently large $n$ \cite{Kuhne2023}. Its Betti numbers $b_m(\mathcal{A}_n)$---the coefficients of the characteristic polynomial of the underlying matroid---are therefore not merely invariants of one example, but encode information about rational arrangements in general.
	
	A second, seemingly unrelated notion is that of a \emph{balanced collection}. A family $\Phi\subseteq 2^{[n]}$ is balanced if there exist positive weights $(\lambda_S)_{S\in\Phi}$ such that
	\[
	\sum_{S\in\Phi}\lambda_S\mathbf{1}_S=\mathbf{1},
	\]
	where $\mathbf{1}_S$ is the characteristic vector of $S$ and $\mathbf{1}=(1,\dots,1)^{\mathsf T}$. A balanced collection is \emph{minimal} if it contains no proper balanced subcollection. Balanced collections were introduced by Bondareva \cite{Bondareva1963} and Shapley \cite{Shapley1967} in the theory of cooperative games, and minimal balanced collections generalise set partitions: every partition of $[n]$ is a minimal balanced collection.
	
	The link between the two notions is classical: by a theorem of Zaslavsky \cite{Zaslavsky1975}, the number of chambers of a real arrangement equals the sum of its Betti numbers, and for $\mathcal{A}_n$ the chambers are in bijection with the \emph{maximal unbalanced} collections of Billera--Moore \cite{BilleraMoore2011}. More directly, the broken-circuit complex of $\mathcal{A}_n$ can be described in terms of minimal balanced collections; this is the starting point of the present work.
	
	K\"uhne \cite{Kuhne2023} made a decisive structural step. He showed that $b_m(\mathcal{A}_n)$ is a fixed linear combination of Stirling numbers of the second kind,
	\[
	b_m(\mathcal{A}_n)=\sum_{k} c_{m,k}\,S(n+1,k),
	\]
	where the coefficients $c_{m,k}$ depend only on $m$ and $k$, not on $n$. The proof introduces \emph{prototypes}: abstract combinatorial objects that encode the local structure of circuits and broken circuits, independently of $n$. K\"uhne further subdivided prototypes into \emph{functional} and \emph{broken}, and reduced the computation of $b_m(\mathcal{A}_n)$ to counting functional prototypes.
	
	However, \cite{Kuhne2023} leaves the counting of functional prototypes as a black box: the prototypes are defined through partitions and maps, and no practical method is given for estimating their number. This is precisely the gap that the present paper fills.
	
	Our approach is based on the combinatorics of minimal balanced collections, developed in our previous work \cite{BludovZuev2025}. There we studied the matrix representation of minimal balanced collections: an $n\times m$ $0$--$1$ matrix $M$ belongs to a certain class $\mathcal{I}_{n,m}^{>0}$ precisely when its columns are distinct and it admits a positive weight vector $\lambda$ with $M\lambda=\mathbf{1}$. Such matrices are in bijection with ordered minimal balanced collections. In the present paper we use this matrix perspective to give a new, equivalent definition of prototypes (Definition~\ref{def:matrix-prototype}), which is both closer to the geometry of $\mathcal{A}_n$ and amenable to counting.
	
	Our first main result is the asymptotic of $b_m(\mathcal{A}_n)$ for every fixed $m$ and variable $n$:
	\[
	b_m(\mathcal{A}_n)=\frac{1}{m!}2^{nm}-O\bigl((3\cdot 2^{m-2})^n\bigr)
	\qquad(n\to\infty).
	\]
	Our second main result refines the correction term: we show that the exponential order above is attained and determine its coefficient,
	\[
	b_m(\mathcal{A}_n)
	=
	\frac{2^{nm}}{m!}
	-
	\frac{m+2}{8(m-2)!}
	\bigl(3\cdot2^{m-2}\bigr)^n
	+
	O_m\!\left(
	\bigl(3\cdot2^{m-2}-1\bigr)^n
	\right)
	\qquad(n\to\infty).
	\]
	The proof combines K\"uhne's theory of prototypes with the matrix-based description of minimal balanced collections, and with a classification of the extremal prototypes that attain the bound $3\cdot2^{m-2}$.
	
	The paper is organised as follows. Section~\ref{sec:prelim} fixes notation and recalls the necessary background on matroids, Betti numbers, and balanced collections. Section~\ref{sec:minimal} develops the matrix theory of minimal balanced collections and weight vectors. Section~\ref{sec:prototypes} introduces prototypes in the matrix language and proves the key bound on broken prototypes. Section~\ref{sec:asymptotic} derives the first asymptotic formula. Section~\ref{sec:constant} refines it, identifying the extremal prototypes and computing the coefficient of the correction term.
	
	\section{Preliminaries}\label{sec:prelim}
	
	\subsection{Notation}
	
	For a positive integer $d$, write $[d]=\{1,\dots,d\}$ and let $\mathbf{1}_d=(1,\dots,1)^{\mathsf T}\in\mathbb{R}^d$; we abbreviate $\mathbf{1}_n$ to $\mathbf{1}$. For $S\subseteq[d]$, let $\mathbf{1}_S\in\{0,1\}^d$ be its characteristic vector. We identify a subset of $[d]$ with its characteristic vector whenever convenient.
	
	Let $\mathcal{M}_{n,m}$ be the set of all $n\times m$ matrices with entries in $\{0,1\}$, viewed as real matrices. For $M\in\mathcal{M}_{n,m}$, its columns are denoted $M_{\ast 1},\dots,M_{\ast m}$ and its rows $M_{1\ast},\dots,M_{n\ast}$.
	
	\subsection{Matroids, Betti numbers, and the NBC complex}
	
	We recall only the minimal background needed; see \cite{Oxley2011} for a thorough treatment.
	
	\begin{definition}\label{def:matroid}
		A \emph{matroid} $M$ on a finite ground set $B$ is a pair $(B,\mathcal{I})$ where $\mathcal{I}\subseteq 2^B$ satisfies the usual independence axioms. The \emph{rank} $\rk(V)$ of $V\subseteq B$ is the maximum size of an independent subset of $V$. A \emph{circuit} is a minimal dependent set.
	\end{definition}
	
	In this paper all matroids arise from arrangements: if $B$ is a finite set of vectors in $\mathbb{R}^n$, the independent sets are the linearly independent subsets of $B$. The \emph{resonance arrangement} $\mathcal{A}_n$ has as its ground set the nonempty $0$--$1$ vectors in $\mathbb{R}^n$, and its matroid is the one determined by linear dependence over $\mathbb{R}$.
	
	\begin{definition}\label{def:nbc}
		Fix an ordering $<$ of $B$. For a circuit $C$, the set $C\setminus\{\max_< C\}$ is a \emph{broken circuit}. A subset of $B$ is \emph{no-broken-circuit} (NBC) if it contains no broken circuit. The family of NBC sets is a simplicial complex, the \emph{NBC complex} $\mathrm{NBC}(M)$.
	\end{definition}
	
	\begin{definition}\label{def:betti}
		The \emph{characteristic polynomial} of a matroid $M$ on $B$ is
		\[
		\chi_M(t)=\sum_{V\subseteq B}(-1)^{|V|}t^{\rk(B)-\rk(V)}.
		\]
		Write
		\[
		\chi_M(t)=\sum_{m=0}^{\rk(B)}(-1)^m b_m(M)t^{\rk(B)-m}.
		\]
		The integers $b_m(M)$ are the \emph{Betti numbers} of $M$.
	\end{definition}
	
	\begin{theorem}[\cite{Zaslavsky1975}]\label{thm:nbc-betti}
		For any matroid $M$,
		\[
		\chi_M(t)=\sum_{V\in\mathrm{NBC}(M)}(-1)^{|V|}t^{\rk(M)-|V|}.
		\]
		Consequently, $b_m(M)=f_{m-1}(\mathrm{NBC}(M))$, where $f_{k}$ is the number of $k$-dimensional faces.
	\end{theorem}
	
	\subsection{Balanced collections and their matrix representation}
	
	\begin{definition}\label{def:balanced}
		A family $\Phi\subseteq 2^{[n]}$ is \emph{balanced} if there exist positive weights $(\lambda_S)_{S\in\Phi}$ such that
		\[
		\sum_{S\in\Phi}\lambda_S\mathbf{1}_S=\mathbf{1}.
		\]
		It is \emph{minimal balanced} if no proper subcollection is balanced.
	\end{definition}
	
	Let $B_{n,m}$ be the number of minimal balanced collections on $[n]$ consisting of exactly $m$ subsets, and let $B_n=\sum_m B_{n,m}$.
	
	For $M\in\mathcal{M}_{n,m}$ of rank $m$, the equation $M\lambda=\mathbf{1}$ has at most one solution. Define
	\[
	\mathcal{I}_{n,m}
	:=\{M\in\mathcal{M}_{n,m}:\rk M=m
	\text{ and }M\lambda=\mathbf{1}\text{ for some }\lambda\in\mathbb{R}^m\}.
	\]
	Consider $m \le n$. For $M\in\mathcal{I}_{n,m}$ we write $\lambda(M)$ for the unique solution and call it the \emph{weight vector}. Put
	\[
	\mathcal{I}_{n,m}^{\neq0}
	:=\{M\in\mathcal{I}_{n,m}:\lambda(M)_j\neq0\ \forall j\},
	\qquad
	\mathcal{I}_{n,m}^{>0}
	:=\{M\in\mathcal{I}_{n,m}:\lambda(M)_j>0\ \forall j\}.
	\]
	
	\begin{proposition}\label{prop:matrix-correspondence}
		If $M\in\mathcal{I}_{n,m}^{>0}$, then its columns are distinct and $\Phi(M)=\{S_1,\dots,S_m\}$ (where $M_{\ast j}=\mathbf{1}_{S_j}$) is a minimal balanced collection. Conversely, every minimal balanced collection $\Phi=\{S_1,\dots,S_m\}$ arises in this way from any ordering of its members.
	\end{proposition}
	
	\begin{corollary}\label{cor:Bnm}
		$B_{n,m}=|\mathcal{I}_{n,m}^{>0}|/m!$.
	\end{corollary}
	
	\section{Minimal balanced collections and weight vectors}\label{sec:minimal}
	
	In this section we develop the matrix theory of minimal balanced collections, following \cite{BludovZuev2025}. Our aim is to isolate the set of \emph{positive weight vectors} $\Lambda_m$, which will be the combinatorial input for the prototype bounds in Section~\ref{sec:prototypes}. We also record the link between minimal balanced collections and circuits of $\mathcal{A}_n$, which serves as the geometric motivation for the framework.
	
	\subsection{Weight vectors and unifying rows}
	
	For $\lambda\in\mathbb{R}^m$, let
	\[
	U(\lambda):=\{u\in\{0,1\}^m:u\lambda=1\}.
	\]
	We call the elements of $U(\lambda)$ the \emph{unifying rows} of $\lambda$. If $A\in\mathcal{I}_{n,m}$ has weight vector $\lambda$, then every row of $A$ belongs to $U(\lambda)$.
	
	\begin{lemma}\label{lem:U-rank}
		Let $\lambda\in\mathbb{R}^m$ and $n\geq m$. There exists $A\in\mathcal{I}_{n,m}$ with $\lambda(A)=\lambda$ if and only if $\rk U(\lambda)=m$. If, moreover, all coordinates of $\lambda$ are nonzero (respectively positive), then such $A$ lies in $\mathcal{I}_{n,m}^{\neq0}$ (respectively $\mathcal{I}_{n,m}^{>0}$).
	\end{lemma}
	
	\begin{proof}
		If $A\in\mathcal{I}_{n,m}$ has weight vector $\lambda$, then every row of $A$ is in $U(\lambda)$, so
		\[
		m=\rk A\leq\rk U(\lambda)\leq m.
		\]
		Conversely, if $\rk U(\lambda)=m$, choose linearly independent $u_1,\dots,u_m\in U(\lambda)$ and, if $n>m$, repeat some of them to obtain $n$ rows. The resulting matrix $A$ has rank $m$ and satisfies $A\lambda=\mathbf{1}$, hence $A\in\mathcal{I}_{n,m}$ and $\lambda(A)=\lambda$. The last assertion is immediate from the definitions.
	\end{proof}
	
	It follows that the set of positive weight vectors is independent of $n\geq m$. We denote it by
	\[
	\Lambda_m:=\{\lambda\in\mathbb{R}_{>0}^m:\rk U(\lambda)=m\}.
	\]
	Equivalently, $\Lambda_m=\{\lambda(A):A\in\mathcal{I}_{n,m}^{>0}\}$ for any $n\geq m$.
	
	For each integer $k$ with $1\leq k<m$, the vector $k^{-1}\mathbf{1}_m$ belongs to $\Lambda_m$: indeed, $U(k^{-1}\mathbf{1}_m)$ consists of all $k$-element subsets of $[m]$, whose characteristic vectors span $\mathbb{R}^m$.
	
	\begin{lemma}\label{lem:minU}
		For every $m\geq1$, $\min_{\lambda\in\Lambda_m}|U(\lambda)|=m$.
	\end{lemma}
	
	\begin{proof}
		Every $\lambda\in\Lambda_m$ satisfies $\rk U(\lambda)=m$, so $|U(\lambda)|\geq m$. Equality is attained at $\lambda=\mathbf{1}_m$, for which $U(\lambda)$ is the set of standard basis vectors.
	\end{proof}
	
	\begin{lemma}\label{lem:maxU}
		For every $m\geq1$,
		\[
		\max_{\lambda\in\Lambda_m}|U(\lambda)|=\binom{m}{\lceil m/2\rceil}.
		\]
	\end{lemma}
	
	\begin{proof}
		Let $r=\lceil m/2\rceil$. For $m\geq2$, the vector $\lambda=r^{-1}\mathbf{1}_m$ lies in $\Lambda_m$, and $U(\lambda)$ is the family of all $r$-element subsets of $[m]$, of size $\binom{m}{r}$. The case $m=1$ is immediate.
		
		For the upper bound, fix $\lambda\in\Lambda_m$. Since all coordinates of $\lambda$ are positive, no two distinct members of $U(\lambda)$ are comparable by inclusion: if $v\subsetneq u$, then $(u-v)\lambda>0$, whereas $u\lambda=v\lambda=1$. Thus $U(\lambda)$ is an antichain in the Boolean lattice, and Sperner's theorem gives $|U(\lambda)|\leq\binom{m}{\lceil m/2\rceil}$.
	\end{proof}
	
	\subsection{Circuits and minimal balanced collections}
	
	\begin{proposition}[\cite{BludovZuev2025}]\label{prop:resonance-circuits}
		For $n\geq1$ and $3\leq m\leq n+1$, let $\kappa_{n,m}$ denote the number of $m$-element circuits of $\mathcal{A}_n$. Then $B_{n+1,m}=2\kappa_{n,m}$. Consequently, $B_{n+1}=2^n+2\kappa_n$, where $\kappa_n=\sum_m\kappa_{n,m}$.
	\end{proposition}
	
	\begin{proof}[Sketch]
		Let $\{\mathbf{1}_{T_1},\dots,\mathbf{1}_{T_m}\}$ be a circuit of $\mathcal{A}_n$. Its linear dependence is unique up to multiplication by a nonzero scalar, so write it as $\sum_{i=1}^m a_i\mathbf{1}_{T_i}=0$. Every $a_i$ is nonzero and both signs occur. Define subsets of $[n+1]$ by
		\[
		S_i:=
		\begin{cases}
			T_i, & a_i>0,\\
			([n]\setminus T_i)\cup\{n+1\}, & a_i<0.
		\end{cases}
		\]
		If $c:=\sum_{a_i<0}|a_i|$, then a coordinate calculation gives
		\[
		\sum_i |a_i|\mathbf{1}_{S_i}=c\mathbf{1}_{n+1},
		\]
		so $\{|a_i|/c\}$ are balancing weights for $\{S_1,\dots,S_m\}$; minimality follows from that of the circuit. Replacing each $a_i$ by $-a_i$ gives the complementary minimal balanced collection $\{[n+1]\setminus S_1,\dots,[n+1]\setminus S_m\}$.
		
		Conversely, let $\Phi=\{S_1,\dots,S_m\}$ be a minimal balanced collection on $[n+1]$, where $m\geq3$, with balancing weights $\lambda_i>0$. Set $T_i:=S_i$ if $n+1\notin S_i$ and $T_i:=[n+1]\setminus S_i$ if $n+1\in S_i$. The $T_i$ are distinct nonempty subsets of $[n]$. For the first $n$ coordinates the balancing equation reads
		\[
		\sum_i \lambda_i\mathbf{1}_{S_i\cap[n]}=\mathbf{1}_n,
		\]
		while for the last coordinate it reads $\sum_{n+1\in S_i}\lambda_i=1$. Substituting $\mathbf{1}_{S_i\cap[n]}=\mathbf{1}_{T_i}$ for $n+1\notin S_i$ and $\mathbf{1}_{S_i\cap[n]}=\mathbf{1}_n-\mathbf{1}_{T_i}$ for $n+1\in S_i$, we obtain
		\[
		\sum_{n+1\notin S_i}\lambda_i\mathbf{1}_{T_i}
		-\sum_{n+1\in S_i}\lambda_i\mathbf{1}_{T_i}=0.
		\]
		Indeed, a proper dependent subset of $\{\mathbf{1}_{T_i}\}$ would yield a proper balanced subcollection of $\Phi$, contradicting its minimality. Since a circuit dependence has exactly two opposite sign patterns, these constructions give a bijection between $m$-element circuits of $\mathcal{A}_n$ and complementary pairs of $m$-element minimal balanced collections on $[n+1]$. Summing over $m$ and using $B_{n+1,1}=1$, $B_{n+1,2}=2^n-1$, and the absence of one- and two-element circuits completes the proof.
	\end{proof}
	
	\section{Prototypes and broken circuits}\label{sec:prototypes}
	
	We now introduce prototypes. K\"uhne's original definition \cite{Kuhne2023} uses partitions and maps; in his convention the zero row is excluded from the prototype and handled separately. In our matrix framework, which is equivalent but more convenient for the present counting arguments, the zero row is counted if it occurs among the distinct rows, and the order of the rows is part of the data.
	
	\subsection{Two equivalent definitions}
	
	Let $A_1,\dots,A_i$ be pairwise distinct nonempty subsets of $[n]$. Taking all intersections of the sets $A_j$ and their complements yields a partition $\pi=\{P_1,\dots,P_k\}$ of $[n+1]$ with $i+1\leq k\leq 2^i$; the block containing $n+1$ is precisely the set of elements of $[n]$ contained in none of the $A_j$. Order the blocks by their binary representation; then $n+1\in P_k$. The tuple $(A_1,\dots,A_i)$ together with $\pi$ determines a map
	\[
	f\colon[k-1]\to\mathcal{P}([i])\setminus\{\varnothing\},
	\qquad
	f(\ell)=\{j\in[i]:P_\ell\subseteq A_j\}.
	\]
	This map is injective because the $A_j$ are distinct, and it is \emph{weakly surjective}: every $j\in[i]$ lies in some $f(\ell)$, since $A_j\neq\varnothing$. Following \cite{Kuhne2023}, an injective weakly surjective map $f\colon[k-1]\to\mathcal{P}([i])\setminus\{\varnothing\}$ is an \emph{$(i,k)$-prototype} in K\"uhne's convention. In our convention, we record the ordered tuple of distinct rows of the matrix with columns $\mathbf{1}_{A_1},\dots,\mathbf{1}_{A_i}$; this tuple may or may not contain the zero row, and its order is inherited from the ordering of the partition blocks.
	
	Equivalently, in the matrix language, let $M\in\mathcal{M}_{n,i}$ have columns $\mathbf{1}_{A_1},\dots,\mathbf{1}_{A_i}$, and let $\mathrm{Row}(M)$ be the ordered tuple of distinct rows of $M$, ordered by the binary representation of the corresponding partition blocks. If the columns of $M$ are distinct and nonzero, then $\mathrm{Row}(M)$ determines an $(i,k)$-prototype with $k=|\mathrm{Row}(M)|$.
	
	\begin{definition}\label{def:matrix-prototype}
		An \emph{$(i,k)$-prototype} is a $k$-tuple $R=(v_1,\dots,v_k)$ of distinct vectors in $\{0,1\}^i$ such that
		\begin{enumerate}
			\item for every coordinate $j\in[i]$ there exists $v\in R$ with $v_j=1$ (equivalently: the columns of the matrix with rows $R$ are all nonzero);
			\item the columns of the $k\times i$ matrix with rows $R$ are pairwise distinct.
		\end{enumerate}
		The vectors $v_1,\dots,v_k$ are the \emph{rows} of the prototype; the order of the tuple is part of the data. The zero vector may or may not belong to $R$.
	\end{definition}
	
	A \emph{realization} of an $(i,k)$-prototype $R=(v_1,\dots,v_k)$ on $[n]$ is a surjective function $\psi\colon[n]\to[k]$; the corresponding $n\times i$ matrix $M$ has rows $M_{x\ast}=v_{\psi(x)}$ for $x\in[n]$. Equivalently, a realization is an $n\times i$ $0$--$1$ matrix whose multiset of rows is $\{v_1,\dots,v_k\}$ with every row occurring at least once, together with the ordering of its distinct rows induced by $R$.
	
	\begin{remark}
		The matrix definition is equivalent to K\"uhne's partition-and-map definition. Namely, put $R^\circ=R\setminus\{0\}$ and let $\varepsilon=1$ if $0\in R$ and $\varepsilon=0$ otherwise. If $R^\circ=(v_1,\dots,v_{k-\varepsilon})$, define
		\[
		f(v_\ell)=\{j\in[i]:(v_\ell)_j=1\}.
		\]
		Then $f$ is an injective weakly surjective map from an ordered set of size $k-\varepsilon$ to $\mathcal{P}([i])\setminus\{\varnothing\}$. Conversely, any such map together with the choice of whether to add the zero row gives a matrix prototype, provided the columns remain distinct and nonzero. Thus the two conventions differ only by the possible presence of the zero row and by the explicit recording of the order, which is canonical in both conventions.
	\end{remark}
	
	\subsection{Circuits and broken circuits}
	
	Fix a total order $<$ on the ground set of $\mathcal{A}_n$; for definiteness, order the nonempty $0$--$1$ vectors lexicographically. The relevance of prototypes is that the property of being a circuit, or of containing a broken circuit, depends only on the prototype. For circuits the order of the rows is irrelevant; for broken circuits it is essential.
	
	\begin{lemma}[\cite{Kuhne2023}]\label{lem:prototype-property}
		Let $(A_1,\dots,A_i)$ be a circuit (respectively, contain a broken circuit) in $\mathcal{A}_n$, and let $P$ be its ordered $(i,k)$-prototype. If $(A'_1,\dots,A'_i)$ in $\mathcal{A}_{n'}$ has the same ordered prototype $P$, then $(A'_1,\dots,A'_i)$ is also a circuit (respectively, contains a broken circuit).
	\end{lemma}
	
	Thus we may speak of \emph{circuit prototypes} and \emph{broken-circuit prototypes}: an ordered $(i,k)$-prototype is a circuit prototype (respectively, a broken-circuit prototype) if some (equivalently, every) realization is a circuit (respectively, contains a broken circuit). For circuit prototypes, the property is invariant under permutation of the rows; for broken-circuit prototypes, it is not.
	
	An $(m,k)$-prototype $R$ \emph{contains} an $(i,l)$-circuit prototype $C$ if there is a set $I\subseteq[m]$ of $i$ coordinates such that the ordered tuple of distinct restrictions $r|_I$, $r\in R$, is exactly $C$. In this case we also say that $R$ is a good expansion of $C$ to length $m$.
	
	\begin{definition}\label{def:containing-broken}
		Let $R$ be an $(m,k)$-prototype. We call $R$ a \emph{containing broken-circuit prototype} if some realization of $R$ contains a broken circuit. By Lemma~\ref{lem:prototype-property}, this property is independent of the choice of realization. Equivalently, $R$ arises by taking a containing circuit $(m+1,\_)$-prototype $Q$, choosing a circuit of $Q$ with $c$ columns, where $3\leq c\leq m+1$, and deleting the column corresponding to the maximal element of that circuit.
	\end{definition}
	
	\begin{remark}
		A good expansion of a broken-circuit prototype need not be a containing broken-circuit prototype. Indeed, being a good expansion is preserved under permutation of the rows, whereas the property of containing a broken circuit is not. Thus the class of containing broken-circuit prototypes is strictly smaller than the class of all good expansions of broken-circuit prototypes.
	\end{remark}
	
	\subsection{Bounding the size of a containing broken-circuit prototype}
	
	The following bound is the technical heart of the asymptotic formula.
	
	\begin{lemma}\label{lem:prototype-bound}
		For every $m>1$, every containing broken-circuit $(m,k)$-prototype has at most $3\cdot 2^{m-2}$ rows; that is, $k\leq 3\cdot 2^{m-2}$.
	\end{lemma}
	
	\begin{proof}
		Let $R$ be a containing broken-circuit $(m,k)$-prototype. By Definition~\ref{def:containing-broken}, $R$ is obtained from a containing circuit $(m+1,\_)$-prototype $Q$ by deleting the maximal circuit-shaping column. Let $c$ be the number of columns of the circuit that is contained in $Q$, so that $3\leq c\leq m+1$, and put $i=c-1$. Thus the broken circuit has $i$ columns, where $2\leq i\leq m$.
		
		Let $R_C$ be the ordered tuple of rows of the circuit prototype $C$ of size $c$, and let $R_B$ be the ordered tuple of rows of the broken circuit obtained by deleting the maximal column of $C$. Deleting a column cannot increase the number of distinct rows, so
		\[
		|R_B|\leq |R_C|.
		\]
		The circuit $C$ corresponds, under column complementation, to a minimal balanced collection with $c$ members. The rows of the corresponding matrix are unifying rows for its positive weight vector. By Lemma~\ref{lem:maxU}, the number of distinct unifying rows is at most $\binom{c}{\lceil c/2\rceil}$. Since column complementation preserves equality of rows, we obtain
		\[
		|R_C|\leq \binom{c}{\lceil c/2\rceil}=\binom{i+1}{\lceil(i+1)/2\rceil}.
		\]
		
		Now reorder the $m$ columns of $R$ so that the $i$ columns of the broken circuit come first. Every row of $R$ restricts to a row of the broken circuit, and for a fixed restriction there are at most $2^{m-i}$ possible suffixes of length $m-i$. Therefore
		\[
		k=|\mathrm{Row}(R)|\leq |R_B|\,2^{m-i}
		\leq \binom{i+1}{\lceil(i+1)/2\rceil}2^{m-i}.
		\]
		Maximising over $2\leq i\leq m$ gives
		\[
		\max_{2\leq i\leq m}\binom{i+1}{\lceil(i+1)/2\rceil}\,2^{m-i}
		=3\cdot 2^{m-2},
		\]
		with equality at $i=2$ and $i=3$. Hence $k\leq 3\cdot 2^{m-2}$.
	\end{proof}
	
	\section{Asymptotic of Betti numbers}\label{sec:asymptotic}
	
	We can now prove the main result.
	
	\begin{theorem}\label{thm:asymptotic}
		For every fixed $m\geq2$,
		\[
		b_m(\mathcal{A}_n)=\frac{1}{m!}2^{nm}-O\bigl((3\cdot 2^{m-2})^n\bigr)
		\qquad(n\to\infty).
		\]
	\end{theorem}
	
	\begin{proof}
		By Theorem~\ref{thm:nbc-betti}, $b_m(\mathcal{A}_n)$ counts the $m$-element no-broken-circuit subsets of the ground set $E_n=\{0,1\}^n\setminus\{0\}$ of $\mathcal{A}_n$. Hence
		\[
		b_m(\mathcal{A}_n)=\binom{2^n-1}{m}-D_m(n),
		\]
		where $D_m(n)$ is the number of $m$-element subsets of $E_n$ that contain a broken circuit.
		
		Adapting K\"uhne's prototype decomposition \cite{Kuhne2023} to containing broken circuits, we obtain
		\[
		D_m(n)=\sum_{k=m}^{2^m}\frac{\#\text{containing broken circuits $(m, k)$-prototypes}}{m!}\,S(n,k),
		\]
		where $S(n,k)$ is a Stirling number of the second kind. The sum is finite because, for fixed $m$ and $k$, there are only finitely many $(m,k)$-prototypes (each is a $k$-tuple of distinct vectors in $\{0,1\}^m$), and their number does not depend on $n$.
		
		By Lemma~\ref{lem:prototype-bound}, every broken $(m,k)$-prototype satisfies $k\le k_0:=3\cdot2^{m-2}$. Therefore the sum over $k$ truncates at $k_0$:
		\[
		D_m(n)=\sum_{k=m}^{k_0}\frac{\#\text{containing broken circuits $(m, k)$-prototypes}}{m!}\,S(n,k).
		\]
		For each $k\le k_0$ we have $S(n,k)=O(k^n)$. Consequently,
		\[
		D_m(n)=O\bigl(k_0^n\bigr)=O\bigl((3\cdot2^{m-2})^n\bigr).
		\]
		On the other hand,
		\[
		\binom{2^n-1}{m}=\frac{2^{nm}}{m!}+O\bigl(2^{n(m-1)}\bigr).
		\]
		Since $3\cdot2^{m-2}\ge2^{m-1}$ for $m\ge2$, the error term $O(2^{n(m-1)})$ is absorbed by $O((3\cdot2^{m-2})^n)$. Combining the estimates gives
		\[
		b_m(\mathcal{A}_n)
		=\frac{2^{nm}}{m!}+O\bigl(2^{n(m-1)}\bigr)
		-O\bigl((3\cdot2^{m-2})^n\bigr)
		=\frac{2^{nm}}{m!}-O\bigl((3\cdot2^{m-2})^n\bigr).
		\]
		This completes the proof.
	\end{proof}
	
	\section{Clarifying the second term}\label{sec:constant}
	
	Theorem~\ref{thm:asymptotic} bounds the correction to
	$2^{nm}/m!$ by a term of order
	$(3\cdot 2^{m-2})^n$. In this section we show that this
	exponential order is attained and determine its coefficient.
	
	Throughout this section, put
	\[
	E_n=\{0,1\}^n\setminus\{0\}.
	\]
	Thus the subsets counted by $b_m(\mathcal A_n)$ are
	$m$-element subsets of $E_n$.
	
	For a prototype $P$, write $\underline P$ for its set of
	rows, forgetting their order. When considering circuit
	prototypes up to isomorphism, we allow permutations of
	both rows and columns. This convention is used only at
	the circuit stage; the maximal element of a realized
	circuit is always taken with respect to the fixed order
	on $E_n$.
	
	If $P$ has $c$ columns and $k$ distinct rows, define
	\[
	\operatorname{Aut}(P)
	=
	\{(\sigma,\tau)\in S_c\times S_k:
	\tau\cdot P\cdot\sigma^{-1}=P\}.
	\]
	Since the rows are distinct, a column permutation
	preserving $\underline P$ determines the corresponding
	row permutation uniquely. Hence $\operatorname{Aut}(P)$
	can equivalently be identified with the stabilizer of
	$\underline P$ under column permutations.
	
	\begin{lemma}\label{lem:realizations}
		Let $\mathcal{Q}$ be an equivalence class of circuit $(i + 1,k)$-prototypes under permutations of columns, and let $\mathcal{P}$ be the corresponding equivalence class of broken circuit $(i,k)$-prototypes under permutations of columns. Let $m \ge i$, and put $r=m-i$ and $q=k2^r$. Then the number of $m$-element subsets
		$S\subseteq E_n$ containing a broken circuit corresponding to some $P \in \mathcal{P}$ is
		\[
		\frac{q!\,S(n,q)}{|\operatorname{Aut}(Q)|\,r!}
		=
		\frac{q^n}{|\operatorname{Aut}(Q)|\,r!}
		+
		O_{m}\bigl((q-1)^n\bigr),
		\]
		where $S(n,q)$ is a Stirling number of the second kind and $Q \in \mathcal{Q}$.
	\end{lemma}
	
	\begin{proof}
		Consider arbitrary $P \in \mathcal{P}$
		
		\textbf{Step 1 (adding the maximal column).}
		Let $Q\in\mathcal{Q}$ be the circuit $(i+1,k)$-prototype from
		which $P$ is obtained by deleting the maximal circuit column.
		In the reverse direction we add one column, which accounts for the factor $i+1$.
		
		\textbf{Step 2 (extending the prototype).}
		Extend $Q$ by $m + 1 - (i + 1) = r$ new columns. Since the $k$ rows of $Q$
		are distinct, each row of $Q$ admits $2^r$ distinct suffixes on
		the new columns, and distinct rows of $Q$ give distinct extended
		rows because they already differ on the original $i + 1$ columns.
		Hence the full expansion of $Q$ has exactly $k\cdot 2^r=q$
		distinct rows. Moreover, every column of the full expansion is
		nonzero and the columns remain pairwise distinct.
		
		\textbf{Step 3 (positioning the circuit columns).}
		Choose the $i+1$ positions of the circuit columns among the
		$m+1$ columns of the expanded prototype; this can be done in
		$\binom{m+1}{i+1}$ ways.
		
		\textbf{Step 4 (ordering the rows).}
		The $q$ distinct rows of the full expansion may be permuted
		arbitrarily, giving a factor $q!$.
		
		\textbf{Step 5 (deleting the maximal column).} To pass from the circuit to a broken circuit, delete the column corresponding to the maximal element of the realized circuit; this accounts for a factor $\frac{1}{m+1}$.
		
		\textbf{Step 6 (realizations).}
		A realization of the resulting $(m,q)$-prototype on $[n]$ is a
		surjective map $\psi\colon[n]\to[q]$, and there are $S(n,q)$
		such maps. For fixed $q$,
		\[
		S(n,q)=\frac{q^n}{q!}+O\bigl((q-1)^n\bigr).
		\]
		
		\textbf{Step 7 (forgetting the column numbering).}
		Since we count unordered $m$-element subsets $S$, we divide by
		$m!$.
		
		\textbf{Step 8 (permuting columns).}
		Since we consider $\mathcal{P}$, not only $P$, we multiply by $i!$.
		
		\textbf{Step 9 (automorphisms).}
		Finally, each equivalence class of prototypes under the action
		of $\operatorname{Aut}(Q)$ has been counted
		$|\operatorname{Aut}(Q)|$ times, so we divide by
		$|\operatorname{Aut}(Q)|$.
		
		Combining the factors and using $m=i+r$, we obtain
		\[
		(i + 1)\cdot\binom{m+1}{i+1}\cdot q!\cdot \frac{1}{m + 1} \cdot S(n,q) \cdot \frac{1}{m!}\cdot i!\cdot\frac{1}{|\operatorname{Aut}(Q)|}
		=
		\frac{q!\,S(n,q)}{|\operatorname{Aut}(Q)|\,r!}.
		\]
	\end{proof}
	
	\begin{remark}
		The factor $r!$ accounts for permutations of the
		additional columns, whereas $|\operatorname{Aut}(Q)|$
		accounts for symmetries of the circuit prototype $Q$.
	\end{remark}
	
	Put
	\[
	q_m=3\cdot2^{m-2}.
	\]
	We first identify the prototypes that can attain this
	number of rows.
	
	\begin{lemma}\label{lem:extremal-prototypes}
		Up to permutations of columns and rows, there is exactly
		one circuit $(3,3)$-prototype, namely
		\[
		P_3=\{000,101,011\},
		\]
		and exactly one circuit $(4,6)$-prototype, namely
		\[
		P_4=\{0000,1111,1100,1010,0101,0011\}.
		\]
	\end{lemma}
	
	Note that the equality case $i = 2$ in Lemma~\ref{lem:prototype-bound} corresponds to \(P_3\), and $i=3$ corresponds to \(P_4\).
	
	\begin{proof}
		We use the correspondence between circuits and minimal
		balanced collections established in Section~\ref{sec:minimal}.
		
		For three columns, an antichain of size three in
		$2^{[3]}$ is either
		\[
		\{100,010,001\}
		\quad\text{or}\quad
		\{110,101,011\}.
		\]
		These correspond to the weight vectors
		$\mathbf1_3$ and $\frac12\mathbf1_3$, respectively.
		In either case, complementing the coordinates of the
		row corresponding to the additional coordinate gives,
		up to column permutation,
		\[
		\{000,101,011\}.
		\]
		
		For four columns, an antichain of size six in
		$2^{[4]}$ is necessarily the complete middle level:
		\[
		U(\lambda)
		=
		\{1100,1010,1001,0110,0101,0011\}.
		\]
		The equations defining these unifying rows imply
		$\lambda=\frac12\mathbf1_4$.
		
		The row corresponding to the additional coordinate is
		one of these six rows. All such choices are equivalent
		under permutations of the four columns. Complementing,
		for example, the coordinates of $1001$ gives
		\[
		\{0000,1111,1100,1010,0101,0011\}.
		\]
		This proves both uniqueness statements.
	\end{proof}
	
	\begin{lemma}\label{lem:automorphisms}
		The extremal prototypes satisfy
		\[
		|\operatorname{Aut}(P_3)|=2,
		\qquad
		|\operatorname{Aut}(P_4)|=8.
		\]
	\end{lemma}
	
	\begin{proof}
		For $P_3$, the first two columns can be exchanged,
		whereas the third column is distinguished by its
		column sum. Thus the only column automorphisms are
		the identity and $(1\,2)$.
		
		For $P_4$, the four rows of weight two are the edges
		of the cycle
		\[
		1-2-4-3-1.
		\]
		Hence its column automorphism group is the automorphism
		group of this cycle, which has order eight.
		Equivalently, it is the group preserving the unordered
		partition
		\[
		\bigl\{\{1,4\},\{2,3\}\bigr\}.
		\]
		It is generated by
		\[
		(1\,4),\qquad (2\,3),\qquad (1\,2)(3\,4).
		\]
		In both cases, a column automorphism determines the
		row permutation uniquely.
	\end{proof}
	
	\begin{definition}\label{def:joint-expansion}
		Let $R$ be a prototype. We call $R$ a \emph{joint expansion}
		of the extremal prototypes $P$ and $Q$ if $R$ is an
		expansion of $P$ and also an expansion of $Q$.
	\end{definition}
	
	\begin{lemma}\label{lem:negligible-overlaps}
		For fixed $m\ge2$, every joint expansion of size $m$ of the broken versions $Q_3$ and $Q_4$ of the extremal prototypes $P_3$ and $P_4$ has strictly fewer than $q_m=3\cdot2^{m-2}$ rows.
	\end{lemma}
	
	\begin{proof}
		Let $R$ be a joint expansion of size \(m\) of some broken versions \(Q_3\) and \(Q_4\) of \(P_3\) and \(P_4\), respectively.
		We first claim that $R$ contains no row all of whose entries
		are $1$. Indeed, the restriction of $R$ to the two columns
		that constitute the extremal $(2,3)$-prototype $Q_3$ coincides with $Q_3$ up to the row order, whose rows are $\{00,10,01\}$; in
		particular, the row $11$ does not occur in this restriction.
		Hence a row of $R$ that were identically $1$ on all $m$
		columns would restrict to the forbidden row $11$ on these
		two columns, a contradiction. Thus $R$ has no all-ones row.
		
		Now suppose that $R$ had $q_m$ rows. By
		Lemma~\ref{lem:prototype-bound}, the bound $q_m$ is attained
		only by full expansions of the broken versions of the extremal prototypes, and a
		full expansion of $Q_4$ is obtained from $Q_4$ by adjoining
		$m-3$ additional columns, with all rows of $Q_4$ extended by
		all possible suffixes. By Step~2 of
		Lemma~\ref{lem:realizations}, such a full expansion has
		exactly $6\cdot2^{m-3}=q_m$ distinct rows, and among them is
		the row obtained by extending the all-ones row \(111\), which is present in every broken version of \(P_4\) (since \(P_4\) contains the row \(1111\)), by the all-ones suffix. This row is identically $1$ on
		all $m$ columns, contradicting the previous paragraph.
		Therefore $R$ is not a full expansion of $Q_4$, and since
		$q_m$ rows can only be attained by a full expansion, $R$ has
		strictly fewer than $q_m$ rows, i.e.\
		$|\mathrm{Row}(R)|\le q_m-1$.
	\end{proof}
	
	\begin{theorem}\label{thm:constant}
		For every fixed $m\ge2$,
		\[
		b_m(\mathcal A_n)
		=
		\frac{2^{nm}}{m!}
		-
		\frac{m+2}{8(m-2)!}
		\bigl(3\cdot2^{m-2}\bigr)^n
		+
		O_m\!\left(
		\bigl(3\cdot2^{m-2}-1\bigr)^n
		\right)
		\qquad(n\to\infty).
		\]
	\end{theorem}
	
	\begin{proof}
		Let $D_m(n)$ denote the number of $m$-element subsets of
		$E_n$ containing a broken circuit. By
		Lemma~\ref{lem:prototype-bound}, every containing
		broken-circuit $(m,k)$-prototype has $k\le q_m$, where
		$q_m=3\cdot2^{m-2}$. By
		Lemma~\ref{lem:extremal-prototypes}, the bound $k=q_m$ is
		attained, up to permutations of rows and columns, only by
		the two extremal prototypes $P_3$ and $P_4$; all remaining
		containing broken-circuit prototypes have at most $q_m-1$
		rows, and by Lemma~\ref{lem:realizations} their total
		contribution is $O_m((q_m-1)^n)$.
		
		It remains to compute the contribution of $P_3$ and $P_4$
		and to justify that the two contributions may be added
		without overcounting. For $P_3$, the number of additional
		columns is $m-2$, and Lemma~\ref{lem:realizations} together
		with $|\operatorname{Aut}(P_3)|=2$
		(Lemma~\ref{lem:automorphisms}) gives the marked
		contribution
		\[
		\frac{q_m^n}{2(m-2)!}
		+
		O_m((q_m-1)^n).
		\]
		For $m\ge3$, the prototype $P_4$ has $m-3$ additional
		columns and, with $|\operatorname{Aut}(P_4)|=8$, contributes
		\[
		\frac{q_m^n}{8(m-3)!}
		+
		O_m((q_m-1)^n).
		\]
		For $m=2$, only $P_3$ contributes, since $Q_4$ has \(3\) columns and \(m-3<0\); the same coefficient formula gives $1/2$.
		
		The only remaining point is that the two contributions
		overlap only in $O_m((q_m-1)^n)$ subsets. By
		Lemma~\ref{lem:negligible-overlaps}, every joint expansion of the broken versions of
		$P_3$ and $P_4$ has strictly fewer than $q_m$ rows. Hence, by
		Lemma~\ref{lem:realizations}, the realizations of such joint
		expansions are negligible, and these are precisely the subsets
		counted twice.
		Hence, for $m\ge3$,
		\[
		\begin{aligned}
			D_m(n)
			&=
			\left(
			\frac{1}{2(m-2)!}
			+
			\frac{1}{8(m-3)!}
			\right)q_m^n
			+
			O_m((q_m-1)^n)\\
			&=
			\frac{m+2}{8(m-2)!}\,q_m^n
			+
			O_m((q_m-1)^n),
		\end{aligned}
		\]
		and the same formula holds for $m=2$.
		
		Finally, the NBC interpretation gives
		\[
		b_m(\mathcal A_n)
		=
		\binom{2^n-1}{m}-D_m(n).
		\]
		Since
		\[
		\binom{2^n-1}{m}
		=
		\frac{2^{nm}}{m!}
		+
		O_m\bigl((2^{m-1})^n\bigr)
		\]
		and
		\[
		2^{m-1}\le q_m-1
		\qquad(m\ge2),
		\]
		this additional error is absorbed into
		$O_m((q_m-1)^n)$. This completes the proof.
	\end{proof}
	
	\section*{Declaration of generative AI and AI-assisted technologies in the manuscript preparation process}
	
	During the preparation of this work the authors used an AI-assisted tool in order to polish the language of the manuscript. After using this tool, the authors reviewed and edited the content as needed and take full responsibility for the content of the published article.
	
	\bibliographystyle{unsrtnat}
	\bibliography{references}
	
\end{document}